\documentclass[preprint,review,10pt]{elsarticle}
\usepackage{amsfonts}
\usepackage{mathtools}
\usepackage{amsmath,amssymb}
\usepackage{graphicx}
\usepackage{setspace}
\usepackage[left=2cm,top=2cm,right=2cm]{geometry}
\usepackage{caption}
\usepackage{subcaption}
\usepackage{amsthm}
\newtheorem{corollary}{Corollary}

\newtheorem{definition}{Definition}

\newtheorem{lemma}{Lemma}

\newtheorem{remark}{Remark}

\newtheorem{theorem}{Theorem}
\numberwithin{equation}{section}

\begin{document}
	\title{Convergence of Sz\'{a}sz--Mirakyan--Durrmeyer operators having Laguerre--type weight}
	\author{Prashantkumar G. Patel}
	\ead{prashant225@spuvvn.edu; prashant225@gmail.com}
	\address{Department of Mathematics, Sardar Patel University, Vallabh Vidyanagar-388 120(Gujarat),
		India}
	\fntext[label*]{Corresponding authors}
	
	\begin{abstract}
In this paper, we introduce a new family of Sz\'{a}sz–Mirakyan–Durrmeyer operators defined on the half-line $[0,\infty)$, constructed using Laguerre-type kernels. The algebraic structure and analytical properties of these operators are thoroughly investigated. Explicit closed-form expressions for the moments are derived, along with a differential recurrence relation connecting successive moments. Quantitative estimates on compact intervals are obtained, and weighted approximation results are provided for unbounded functions. Furthermore, the asymptotic behavior of the central moments is analyzed. We establish both local and global $L_p$-convergence results and identify the eigenfunctions associated with these operators. These findings demonstrate the effectiveness of the proposed generalized operators in extending classical approximation results to the unbounded domain.
	\end{abstract}
	\maketitle
	\section{Introduction}
	The problem of approximating continuous functions on unbounded intervals, particularly $[0,\infty)$, has been a central topic in approximation theory since the early 20th century. Classical results such as the Weierstrass approximation theorem guarantee uniform approximation on compact intervals, but these techniques do not directly extend to infinite domains.
	
	To address this, Sz\'{a}sz (1950–1951) \cite{szasz1950generalization} introduced a class of operators, now known as the Sz\'{a}sz operators, defined for a function $f\in C[0,\infty)$ by 
	\begin{eqnarray}\label{SMoperators}
		S_n(f,x):=e^{-nx} \sum_{k=0}^{\infty} \dfrac{(nx)^k}{k!}f\left(\dfrac{k}{n}\right), x\geq 0
	\end{eqnarray}
	Independently, Mirakyan (1941) \cite{mirakjan1941approximation} studied similar constructions using positive linear operators to approximate continuous functions on  $[0,\infty)$. These operators extended the classical Bernstein polynomials, which are effective on bounded intervals $[0,1]$, to the unbounded setting by exploiting Poisson-type weights.\\
	The Sz\'asz–Mirakyan operators are positive, linear, and preserve constants, making them a natural tool for studying approximation on infinite intervals. Over the decades, these operators have been generalized in many directions, including Durrmeyer-type modifications, $q$-analogues, and weighted $L^p$ spaces, providing flexible frameworks for approximating functions with exponential growth or decay at infinity.
	In 1981, Sablonni\`{e}re \cite{sablonniere1981operateurs} introduced an extended version of the Durrmeyer operator by incorporating the classical Jacobi weight function. The operator is defined as
	\begin{eqnarray*}
		B_n^{(\alpha,\beta)} (f,x) = \sum_{k=0}^{\infty} \frac{\langle f, p_{n,k}\rangle^*_{\alpha,\beta}}{\langle 1, p_{n,k} \rangle^*_{\alpha,\beta}} \, p_{n,k} (x),
	\end{eqnarray*}
	where the weighted inner product is
	\[
	\langle f, g\rangle^*_{\alpha,\beta} = \int_0^1 f(t) g(t) t^{\alpha} (1-t)^{\beta} dt, \quad \alpha, \beta > -1, \textrm{ and } p_{n,k}(x) ={n \choose k} x^k (1-x)^{n-k}.
	\]
	For the special case $\alpha = \beta = 0$, the classical Durrmeyer operator is recovered. Recently, Bernstein–Jacobi-type operators preserving derivatives were studied and established by Lara-Velasco and P\'{e}rez \cite{lara2024bernstein}, extending these operators to preserve differential properties in addition to function approximation.
	To approximate integrable functions on the infinite interval $[0,\infty)$, Mazhar and Totik \cite{mazhar1985approximation} introduced a Durrmeyer-type variant of the Sz\'{a}sz operators \eqref{SMoperators} as
	\begin{eqnarray*}
		D_n(f,x) := n\sum_{k=0}^{\infty} \psi_{n,k}(x) \int_0^\infty f(t) \, \psi_{n,k}(t) \, dt,
	\end{eqnarray*}
	where $\psi_{n,k}(t)=\frac{(nx)^k e^{-nx}}{k!}$. Since then, many authors have studied these operators and proposed various generalizations \cite{gupta2014convergence,abel2025asymptotic,patel2024certain,patel2025durrmeyer} and many other.  
	
	However, unlike the Bernstein–Jacobi–Durrmeyer operators on $[0,1]$, these operators have not been extensively investigated with classical Jacobi weights on $[0,\infty)$, primarily because the Jacobi weight 
	\[
	w^*_{\alpha,\beta}(t) = t^\alpha (1-t)^\beta
	\] 
	is naturally defined on a bounded interval $[0,1]$. To extend the operators to the unbounded interval, it is more convenient to use a Laguerre-type weight defined by 
	\[
	\omega_{\alpha, \beta} (t) := e^{\beta t} t^\alpha, \quad \alpha > -1, \, \beta \in \mathbb{R}.
	\]  
	This weight ensures that the integrals converge for a wide class of functions on $[0,\infty)$, preserves the positivity of the operators, and allows us to define moments and study approximation properties in weighted $L_p$ spaces.  Now, let us define the following weighted inner product for integrable function $f$ and $g$ as
	$$\langle f,g\rangle_{\alpha,\beta}:=\int_0^{\infty} f(t)g(t) e^{\beta t}t^\alpha dt, \alpha, \beta \in \mathbb{R}.$$

	These considerations, which motivate one to introduce the following Laguerre-type Durrmeyer operators, allow a natural generalization of the Sz\'{a}sz–Mirakyan operators while retaining the desirable properties of positivity, linearity, and normalization: Fix parameters $\alpha>-1$ and $\beta>0$.
	\begin{eqnarray}\label{eq:operator}
		M_n^{(\alpha,\beta)}[f](x) := \sum_{k=0}^{\infty} \frac{\langle f, \psi_{n,k}\rangle_{\alpha,\beta}}{\langle 1, \psi_{n,k} \rangle_{\alpha,\beta}} \, \psi_{n,k} (x),
	\end{eqnarray}
	
	can be written explicitly as
	\[
	M_n^{(\alpha,\beta)}[f](x) = \sum_{k=0}^{\infty} \left( \frac{(n - \beta)^{k + \alpha + 1}}{\Gamma(k + \alpha + 1)} \int_0^{\infty} f(t) t^{k+\alpha} e^{-(n - \beta)t} dt \right) \psi_{n,k}(x),
	\]
	provided that \( n > \beta \) and \( \alpha > -1 \), so that the integrals converge.\\
	For a real number $A\geq 0$, let $E_{A}\left[ 0,\infty \right) $ be the
	class of all functions $f:\left[ 0,\infty \right) \rightarrow \mathbb{R}$
	which satisfy, for some positive constant $K$, the growth condition $%
	\left\vert f\left( t\right) \right\vert \leq K\exp \left( At\right) $, for $%
	t\geq 0$. Let $E\left[ 0,\infty \right) =\bigcup_{A\geq 0}E_{A}\left[
	0,\infty \right) $, the space of all functions of (at most)\ exponential
	growth. The operators \eqref{eq:operator} interpolate the function $f$ at
	the point $x=0$, i.e. $ M_n^{(\alpha,\beta)}[f]\left(
	0\right) =\frac{(n - \beta)^{\alpha + 1}}{\Gamma(\alpha + 1)} \int_0^{\infty} f(t) \, t^{\alpha} e^{-(n - \beta)t} \, dt$. It is obvious that, the operators $M_{n}^{(\alpha,\beta)}$
	are linear and positive. We observe that the operators are well defined on
	the space $E\left[ 0,\infty \right) $ since $f\in E_{A}\left[ 0,\infty
	\right) $, implies that 
	\begin{equation*}
		\left\vert \left( M_{n}^{\alpha,\beta}f\right) \left( x\right)
		\right\vert \leq K \left( \frac{n - \beta}{n - \beta - A} \right)^{\alpha + 1} \exp\left( \frac{n x A}{n - \beta - A} \right).
	\end{equation*}
In the existing literature, no attempt has yet been made to define operators analogous to the Bernstein--Jacobi--Durrmeyer operators on the unbounded interval. Motivated by the properties and success of the Bernstein--Jacobi--Durrmeyer operators, we introduce a new class of operators \eqref{eq:operator}, referred to as the Sz\'{a}sz--Mirakyan--Laguerre--Durrmeyer operators. 
	
	The main objectives of this work are to establish the moments of these operators, to derive their recurrence relations in differential form, and to investigate their global $(L_p)$-convergence in suitable weighted spaces. Such an analysis provides a deeper understanding of their approximation behavior on the infinite interval $[0,\infty)$ and extends the classical theory of Durrmeyer-type operators in a Laguerre-weighted framework.
	
	\section{Moments of the Operators}
	
	In this section, we compute the moments of the operator \( M_n^{(\alpha,\beta)} \) defined by \eqref{eq:operator}.
Consider $f(t)=t^r$, $r\in\{0,1,2,\dots\}$.
	\begin{align*}
		M_n^{(\alpha,\beta)}[t^r](x)
		&=\sum_{k=0}^{\infty}\frac{(n-\beta)^{k+\alpha+1}}{\Gamma(k+\alpha+1)}\cdot\frac{\Gamma(k+\alpha+r+1)}{(n-\beta)^{k+\alpha+r+1}}\,\psi_{n,k}(x)\\
		&=\frac{1}{(n-\beta)^r}\sum_{k=0}^{\infty}\frac{\Gamma(k+\alpha+r+1)}{\Gamma(k+\alpha+1)}\,\psi_{n,k}(x)
		=\frac{1}{(n-\beta)^r}\sum_{k=0}^{\infty}(\alpha+1+k)_r\,\psi_{n,k}(x).
	\end{align*}
	Here $(a)_r=a(a+1)\cdots(a+r-1)$ denotes the rising factorial (Pochhammer symbol).\\
	Using Gamma identities one checks the identity
	\[ (\alpha+1+k)_r = (\alpha+1)_r\,\frac{(\alpha+r+1)_k}{(\alpha+1)_k},\]
	so that the sum can be expressed via the confluent hypergeometric function $\,_1F_1$:
	\begin{align}\label{eq:closed-form}
		M_n^{(\alpha,\beta)}[t^r](x)
		&=\frac{(\alpha+1)_r}{(n-\beta)^r} \, e^{-nx} \, {}_1F_1(\alpha+r+1;\,\alpha+1;\,nx),\qquad r=0,1,2,\dots
	\end{align}


	Recall the $r$-th moment as
	\[\mu_r^{(\alpha,\beta)}(x):=M_n^{(\alpha,\beta)}[t^r](x)=\frac{1}{(n-\beta)^r}\sum_{k=0}^{\infty}(\alpha+1+k)_r\,\psi_{n,k}(x).
	\]
	Using the identity
	\[\frac{d}{dx}\psi_{n,k}(x)=n\big(\psi_{n,k-1}(x)-\psi_{n,k}(x)\big),\qquad\psi_{n,-1}(x):=0,\]
	and the finite-difference relation for rising factorials
	\[(u+1)_r-(u)_r = r\,(u+1)_{r-1},\]
	one obtains 
	\begin{equation}\label{eq:recurrence-shifted}
		\frac{d}{dx}\,\mu_r^{(\alpha,\beta)}(x)
		=\frac{n r}{\,n-\beta\,}\;\mu_{r-1}^{(\alpha+1,\beta)}(x),\qquad r\ge1.
	\end{equation}
	This is the correct differential recurrence: the $(r-1)$-th moment on the right-hand side is taken with the parameter $\alpha$ replaced by $\alpha+1$.

	In particular the normalization (zeroth moment) is
	\[\mu_0^{(\alpha,\beta)}(x)=M_n^{(\alpha,\beta)}[1](x)=\sum_{k\ge0}\psi_{n,k}(x)=1.
	\]
	The sums can be evaluated directly using the Poisson moment identities $\sum_k\psi_{n,k}=1$, $\sum_k k\psi_{n,k}=nx$, $\sum_k k(k-1)\psi_{n,k}=(nx)^2$, $\sum_k k(k-1)(k-2)\psi_{n,k}=(nx)^3$. The first three nontrivial moments are:
	\begin{align*}
		\mu_1^{(\alpha,\beta)}(x)
		&=\frac{1}{n-\beta}\sum_{k\ge0}(\alpha+1+k)\psi_{n,k}(x)
		=\frac{\alpha+1+nx}{n-\beta},\\
		\mu_2^{(\alpha,\beta)}(x)
		&=\frac{1}{(n-\beta)^2}\Big((nx)^2+(2\alpha+4)nx+(\alpha+1)(\alpha+2)\Big),\\
		\mu_3^{(\alpha,\beta)}(x)
		&=\frac{1}{(n-\beta)^3}\Big((nx)^3+3(\alpha+3)(nx)^2+3(3\alpha+5)(nx)+(\alpha+1)(\alpha+2)(\alpha+3)\Big).\\	\mu_4^{(\alpha,\beta)}(x) &= \frac{(nx)^4 + 6(\alpha+4)(nx)^3 + (11\alpha^2 + 58\alpha + 72)(nx)^2}{(n-\beta)^4} \\
		&\quad + \frac{6(\alpha+1)(\alpha+2)(\alpha+4)\, nx + (\alpha+1)(\alpha+2)(\alpha+3)(\alpha+4)}{(n-\beta)^4}.
	\end{align*}
	
The confluent hypergeometric function ${}_1F_1(a;b;z)$ satisfies the standard contiguous relation
(see \cite[Chapter 13]{daalhuis2010confluent})
\[
r\,{}_1F_1(a-1;b;z)+(\alpha+2r+1+z)\,{}_1F_1(a;b;z)-a\,{}_1F_1(a+1;b;z)=0,
\]
with the substitutions
\[
a=\alpha+r+1,\qquad b=\alpha+1,\qquad z=nx.
\]

Recall the closed form for the \(r\)-th moment
\[
\mu_r^{(\alpha,\beta)}(x)=\frac{(\alpha+1)_r}{(n-\beta)^r}e^{-z}\;{}_1F_1(\alpha+r+1;\alpha+1;z),\qquad z=nx.
\]
Set
\[
F_r(z):={}_1F_1(\alpha+r+1;\alpha+1;z),\qquad
C_r:=\frac{(\alpha+1)_r}{(n-\beta)^r}e^{-z},
\]
so that \(\mu_r^{(\alpha,\beta)}(x)=C_rF_r(z)\).

The contiguous relation rewrites as
\[
rF_{r-1}+(\alpha+2r+1+z)F_r-(\alpha+r+1)F_{r+1}=0.
\]
Multiply this identity by \(C_{r-1}\) and use the relations
\[
\frac{C_r}{C_{r-1}}=\frac{\alpha+r}{\,n-\beta\,},\qquad
\frac{C_{r+1}}{C_{r-1}}=\frac{(\alpha+r+1)(\alpha+r)}{(n-\beta)^2},
\]
equivalently
\[
C_{r-1}F_{r-1}=\mu_{r-1}^{(\alpha,\beta)},\qquad
C_{r-1}F_r=\frac{n-\beta}{\alpha+r}\,\mu_r^{(\alpha,\beta)},\qquad
C_{r-1}F_{r+1}=\frac{(n-\beta)^2}{(\alpha+r)(\alpha+r+1)}\,\mu_{r+1}^{(\alpha,\beta)}.
\]

Substituting these into the multiplied contiguous relation and simplifying yields
\[
r\,\mu_{r-1}^{(\alpha,\beta)}+\frac{(n-\beta)(\alpha+2r+1+z)}{\alpha+r}\,\mu_r^{(\alpha,\beta)}
-\frac{(n-\beta)^2}{\alpha+r}\,\mu_{r+1}^{(\alpha,\beta)}=0.
\]
Multiplying through by \(\alpha+r\) and rearranging gives the three-term recurrence
\begin{equation}\label{eq:recurrence-r-final}
	(n-\beta)^2\,\mu_{r+1}^{(\alpha,\beta)}(x)
		= r(\alpha+r)\,\mu_{r-1}^{(\alpha,\beta)}(x)
		+ (n-\beta)\big(\alpha+2r+1+nx\big)\,\mu_r^{(\alpha,\beta)}(x),
		\;\qquad r=0,1,2,\dots
\end{equation}

Using the binomial expansion
\[
(t-x)^r=\sum_{j=0}^r\binom{r}{j}(-x)^{\,r-j}t^j
\]
and the raw moments $\mu_j^{(\alpha,\beta)}(x)=M_n^{(\alpha,\beta)}[t^j](x)$ we have
\[
M_n^{(\alpha,\beta)}[(t-x)^r](x)
=\sum_{j=0}^r\binom{r}{j}(-x)^{\,r-j}\mu_j^{(\alpha,\beta)}(x).
\]
Substituting the explicit $\mu_j^{(\alpha,\beta)}$ given in the text (with $\mu_0=1$) yields the following closed forms.

\begin{eqnarray*}
	M_n^{(\alpha,\beta)}\big[(t-x)\big](x)	&=&\frac{\alpha+1+\beta x}{\,n-\beta\,},\\
	M_n^{(\alpha,\beta)}\big[(t-x)^2\big](x)&=&\frac{(\alpha+1)(\alpha+2)+2x\big(n+\beta(\alpha+1)\big)+\beta^2 x^2}{(n-\beta)^2}\,,\\
	M_n^{(\alpha,\beta)}\big[(t-x)^3\big](x) &=&\dfrac{\beta^{3}x^{3}+\big(3\alpha\beta^{2}+3\beta^{2}+6\beta n\big)x^{2}
		+\big(3\alpha^{2}\beta-3\alpha^{2}n+9\alpha\beta+6\beta+9n\big)x }{(n-\beta)^3}\\
	&& +\dfrac{(\alpha+1)(\alpha+2)(\alpha+3)}{(n-\beta)^3}\,,\\
	M_n^{(\alpha,\beta)}\big[(t-x)^4\big](x) &=&\dfrac{\beta^{4}x^{4}+\big(4\alpha\beta^{3}+2\alpha n^{3}+4\beta^{3}+12\beta^{2}n+8n^{3}\big)x^{3}}{(n-\beta)^4}\\
	&&+\dfrac{\big(6\alpha^{2}\beta^{2}-12\alpha^{2}\beta n+17\alpha^{2}n^{2}+10\alpha\beta^{2}-40\alpha\beta n+40\alpha n^{2}+12\beta^{2}+36\beta n+24n^{2}\big)x^{2}}{(n-\beta)^4}\\
	&&+\dfrac{\big(4\alpha^{3}\beta+2\alpha^{3}n+24\alpha^{2}\beta+18\alpha^{2}n+44\alpha\beta+40\alpha n+24\beta+24n\big)x}{(n-\beta)^4}\\
	&&+\dfrac{(\alpha+1)(\alpha+2)(\alpha+3)(\alpha+4)}{(n-\beta)^4}\,.
\end{eqnarray*}

For the operator \(M_n^{(\alpha,\beta)}\) recall the raw moments
\[
\mu_j^{(\alpha,\beta)}(x):=M_n^{(\alpha,\beta)}[t^j](x)
=\frac{(\alpha+1)_j}{(n-\beta)^j}\,e^{-nx}\;{}_1F_1(\alpha+j+1;\,\alpha+1;\,nx),
\qquad j=0,1,2,\dots
\]
(where \((\alpha+1)_0:=1\) and \({}_1F_1\) is the confluent hypergeometric function).

Using the binomial expansion \((t-x)^r=\sum_{j=0}^r\binom{r}{j}(-x)^{\,r-j}t^j\), the $r$-th central moment
\[
M_n^{(\alpha,\beta)}\big[(t-x)^r\big](x)
=\sum_{j=0}^r\binom{r}{j}(-x)^{\,r-j}\,\mu_j^{(\alpha,\beta)}(x)
\]
has the closed form
\[
M_n^{(\alpha,\beta)}\big[(t-x)^r\big](x)
	= e^{-nx}\sum_{j=0}^r \binom{r}{j}(-x)^{\,r-j}\frac{(\alpha+1)_j}{(n-\beta)^j}\;
	{}_1F_1(\alpha+j+1;\,\alpha+1;\,nx).%
\]

This formula is valid for every integer \(r\ge 0\).  It expresses each central moment as a finite linear combination (length \(r+1\)) of confluent hypergeometric terms.  

If desired one may factor out \((n-\beta)^{-r}e^{-nx}\) to get
\[
M_n^{(\alpha,\beta)}\big[(t-x)^r\big](x)
=\frac{e^{-nx}}{(n-\beta)^r}\sum_{j=0}^r \binom{r}{j}(-x)^{\,r-j}(\alpha+1)_j\,(n-\beta)^{\,r-j}
\,{}_1F_1(\alpha+j+1;\,\alpha+1;\,nx).
\]
\section{Asymptotic Analysis of Central Moments}

The following result give the asymptotic relation of central moments of the operators $	M_n^{(\alpha,\beta)}$.
\begin{theorem}\label{thm:asymtotic}
	For fixed $x > 0$, $\alpha > -1$, and $\beta \in \mathbb{R}$, the central moments of the operator $M_n^{(\alpha,\beta)}$ satisfy the following asymptotic expansions as $n \to \infty$:
	\begin{enumerate}
		\item For $r = 1$:
		\[
		M_n^{(\alpha,\beta)}[(t-x)](x) = \frac{\alpha + 1 + \beta x}{n}.
		\]
		\item For every integer $r \geq 2$:
		\[
		M_n^{(\alpha,\beta)}[(t-x)^r](x) = \frac{r(r-1)\beta^{r-2}x^{r-1}}{n^{r-1}} + O(n^{-r}).
		\]
	\end{enumerate}
\end{theorem}

\begin{proof}
	The $j$-th moment is given by:
	\[
	\mu_j^{(\alpha,\beta)}(x) = \frac{(\alpha+1)_j}{(n-\beta)^j} e^{-nx} \, {}_1F_1(\alpha+j+1; \alpha+1; nx),
	\]
	where $(\alpha+1)_j$ denotes the Pochhammer symbol.
	
	Introduce the notation:
	\[
	z = nx, \quad A = \frac{n}{n-\beta} = 1 + \frac{\beta}{n-\beta}, \quad \delta = A - 1 = \frac{\beta}{n-\beta}.
	\]
	
	Using the known asymptotic expansion for the confluent hypergeometric function:
	\[
	e^{-z} \, {}_1F_1(a; b; z) \sim \frac{\Gamma(b)}{\Gamma(a)} z^{a-b} \left[1 + \frac{(b-a)(1-a)}{z} + O(z^{-2})\right], \quad z \to +\infty,
	\]
	with $a = \alpha+j+1$, $b = \alpha+1$, we obtain:
	\[
	e^{-z} \, {}_1F_1(\alpha+j+1; \alpha+1; z) \sim \frac{1}{(\alpha+1)_j} z^j \left[1 + \frac{j(\alpha+j)}{z} + O(z^{-2})\right].
	\]
	
	Substituting into the expression for $\mu_j$:
	\begin{align*}
		\mu_j^{(\alpha,\beta)}(x) &= \frac{(\alpha+1)_j}{(n-\beta)^j} \cdot \frac{1}{(\alpha+1)_j} z^j \left[1 + \frac{j(\alpha+j)}{z} + O(z^{-2})\right] \\
		&= \frac{z^j}{(n-\beta)^j} \left[1 + \frac{j(\alpha+j)}{z} + O(z^{-2})\right].
	\end{align*}
	
	Since $z = nx$ and $\frac{z}{n-\beta} = xA$, we arrive at:
	\begin{equation}\label{eq:raw_moment_asymptotic}
		\mu_j^{(\alpha,\beta)}(x) = x^j A^j \left[1 + \frac{j(\alpha+j)}{nx} + O(n^{-2})\right]. 
	\end{equation}
	The $r$-th central moment is:
	\[
	M_n^{(\alpha,\beta)}[(t-x)^r](x) = \sum_{j=0}^r \binom{r}{j} (-x)^{r-j} \mu_j^{(\alpha,\beta)}(x).
	\]
	
	Substituting \eqref{eq:raw_moment_asymptotic}:
	\begin{align*}
		M_n^{(\alpha,\beta)}[(t-x)^r](x) &= \sum_{j=0}^r \binom{r}{j} (-x)^{r-j} x^j A^j \left[1 + \frac{j(\alpha+j)}{nx} + O(n^{-2})\right] \\
		&= x^r \sum_{j=0}^r \binom{r}{j} (-1)^{r-j} A^j \left[1 + \frac{j(\alpha+j)}{nx} + O(n^{-2})\right]. 
	\end{align*}
	
	Define the sums:
	\[
	S_0 = \sum_{j=0}^r \binom{r}{j} (-1)^{r-j} A^j, \quad 
	S_1 = \sum_{j=0}^r \binom{r}{j} (-1)^{r-j} A^j j(\alpha+j).
	\]
	
	Then \eqref{eq:raw_moment_asymptotic} becomes:
	\begin{equation}\label{eq:central_moment_expansion}
		M_n^{(\alpha,\beta)}[(t-x)^r](x) = x^r S_0 + \frac{x^{r-1}}{n} S_1 + O(n^{-r}). 
	\end{equation}
	Consider the generating polynomial:
	\[
	F(z) = \sum_{j=0}^r \binom{r}{j} (-1)^{r-j} A^j z^j = (Az - 1)^r.
	\]
	
	Then	$S_0 = F(1) = (A - 1)^r$ and $S_1 = (\alpha+1)F'(1) + F''(1)$.\\
	Compute derivatives:
	\begin{align*}
		F'(z) &= rA(Az - 1)^{r-1} \Rightarrow F'(1) = rA(A - 1)^{r-1} \\
		F''(z) &= r(r-1)A^2(Az - 1)^{r-2} \Rightarrow F''(1) = r(r-1)A^2(A - 1)^{r-2}
	\end{align*}
	Thus
	\begin{align*}
		S_1 &= (\alpha+1)rA(A - 1)^{r-1} + r(r-1)A^2(A - 1)^{r-2} \\
		&= rA(A - 1)^{r-2} \left[(\alpha+1)(A - 1) + (r-1)A\right] \\
		&= rA(A - 1)^{r-2} \left[(\alpha + r)A - (\alpha + 1)\right]. 
	\end{align*}
	Substituting into \eqref{eq:central_moment_expansion}:
	\begin{equation}\label{eq:two_term}
		M_n^{(\alpha,\beta)}[(t-x)^r](x) = x^r (A - 1)^r  + \frac{x^{r-1}}{n} rA(A - 1)^{r-2} \left[(\alpha + r)A - (\alpha + 1)\right] + O(n^{-r}). 
	\end{equation}
	Now expand $A$ and $(A-1)$ in powers of $n^{-1}$:
	\[
	A = 1 + \frac{\beta}{n} + O(n^{-2}), \quad A - 1 = \frac{\beta}{n} + O(n^{-2}).
	\]
	Then, we note that 
	\begin{eqnarray*}
		(A - 1)^r &=& \frac{\beta^r}{n^r} + O(n^{-r-1})\\
		A(A - 1)^{r-2} &=& \frac{\beta^{r-2}}{n^{r-2}} + O(n^{-r+1})\\
		(\alpha + r)A - (\alpha + 1) &=& (r - 1) + O(n^{-1}).
	\end{eqnarray*}		
		
	Substituting into \eqref{eq:two_term}:
	\begin{align*}
		M_n^{(\alpha,\beta)}[(t-x)^r](x) &= x^r \cdot \frac{\beta^r}{n^r} + \frac{x^{r-1}}{n} \cdot r \cdot \frac{\beta^{r-2}}{n^{r-2}} \cdot (r - 1) + O(n^{-r}) \\
		&= \frac{r(r - 1)\beta^{r-2}x^{r-1}}{n^{r-1}} + O(n^{-r}). \tag{6}
	\end{align*}
	
	For the special case $r = 1$, from \eqref{eq:two_term}:
	\[
	M_n^{(\alpha,\beta)}[(t-x)](x) = (A - 1)x + \frac{1}{n} A[(\alpha + 1)A - (\alpha + 1)] + O(n^{-2}).
	\]
	Since $A - 1 = \frac{\beta}{n} + O(n^{-2})$ and $A = 1 + O(n^{-1})$:
	\[
	M_n^{(\alpha,\beta)}[(t-x)](x) = \frac{\beta x}{n} + \frac{\alpha + 1}{n} + O(n^{-2}) = \frac{\alpha + 1 + \beta x}{n} + O(n^{-2}).
	\]
	
	This completes the proof of Theorem \ref{thm:asymtotic}.

\end{proof}
\section{Quantitative estimate on compact intervals}
In this section, we investigate the rate of convergence of the proposed Szász–Mirakyan–Durrmeyer operators on compact subsets of $[0,\infty)$. Beyond establishing mere convergence, it is essential to obtain quantitative estimates that describe how rapidly the operators approximate a given function. Such results are formulated in terms of the modulus of continuity, which provide sharp bounds on the approximation error. The analysis presented here highlights the efficiency of the operators in achieving uniform approximation on bounded intervals and offers a deeper understanding of their local behavior.\\
Let \([0, a]\) be a compact interval. For \(f \in C([0, a])\), define the usual modulus of continuity:
\[
\omega(f, \delta) = \sup_{\substack{x,t \in [0,a] \\ |t-x| \le \delta}} |f(t)-f(x)|.
\]

\begin{theorem}[Quantitative estimate on compact intervals]
	\label{thm:compact-correct}
	Assume \(n > \beta\) and \(\alpha > -1\). For every \(f \in C([0, a])\) and every \(n > \beta\), there exists a constant \(C > 0\) depending on \(a, \alpha, \beta\) such that
	\[
	\sup_{x \in [0, a]} \left| M_n^{(\alpha,\beta)} f(x) - f(x) \right|
	\le C \, \omega\left(f, \frac{1}{\sqrt{n}}\right).
	\]
	In particular, \(M_n^{(\alpha,\beta)} f \to f\) uniformly on \([0, a]\) as \(n \to \infty\).
\end{theorem}

\begin{proof}
	We use the following inequality for any \(\delta > 0\):
	\[
	|f(t)-f(x)| \le \left(1 + \frac{|t-x|}{\delta}\right) \omega(f, \delta).
	\]
	Applying the operator \(M_n^{(\alpha,\beta)}\) and using its linearity and positivity, we get:
	\[
	|M_n^{(\alpha,\beta)} f(x) - f(x)| \le \omega(f, \delta) \left(1 + \frac{1}{\delta} M_n^{(\alpha,\beta)}(|t-x|; x)\right).
	\]
	By Cauchy-Schwarz, 
	\[
	M_n^{(\alpha,\beta)}(|t-x|; x) \le \sqrt{M_n^{(\alpha,\beta)}((t-x)^2; x)} = \sqrt{\mu_{n,2}(x)}.
	\]
	The correct computation of the second moment is:
	\[
	\mu_{n,2}(x) = \frac{\beta^2 x^2 + [n + 2\beta(\alpha+1)] x + (\alpha+1)(\alpha+2)}{(n-\beta)^2}.
	\]
	For \(x \in [0, a]\), there exists a constant \(C_1 > 0\) depending on \(a, \alpha, \beta\) such that
	\[
	\mu_{n,2}(x) \le \frac{C_1}{n}.
	\]
	Thus,
	\[
	M_n^{(\alpha,\beta)}(|t-x|; x) \le \sqrt{\frac{C_1}{n}}.
	\]
	Choosing \(\delta = \frac{1}{\sqrt{n}}\), we obtain
	\[
	|M_n^{(\alpha,\beta)} f(x) - f(x)| \le \omega\left(f, \frac{1}{\sqrt{n}}\right) \left(1 + \sqrt{C_1}\right).
	\]
	Taking \(C = 1 + \sqrt{C_1}\) completes the proof.
\end{proof}

\section{Weighted approximation for unbounded functions}
To study the approximation behavior of the proposed operators for unbounded functions, we consider the weighted space 
$C_{k,\varphi}([0,\infty))$ is a suitable weight function controlling growth at infinity. In this setting, we establish Korovkin-type convergence results, which ensure that the sequence of Sz\'{a}sz–Mirakyan–Durrmeyer operators preserves approximation in the weighted sense. By verifying the convergence of the operators on a standard test set of functions, we extend the classical Korovkin theorem to the unbounded domain $[0,\infty)$. This framework provides a natural and effective approach for analyzing the weighted approximation of functions exhibiting polynomial growth.\\
Let $\varphi(x)=1+x^2$, and consider the weighted spaces
\[
B_\varphi([0,\infty)) = \{ f:[0,\infty)\to\mathbb{R}: |f(x)|\le M_f \varphi(x) \text{ for some } M_f>0\},
\]
\[
C_\varphi([0,\infty)) = C([0,\infty)) \cap B_\varphi([0,\infty)),
\quad
C_{k,\varphi}([0,\infty)) = \Big\{ f \in C_\varphi([0,\infty)): \lim_{x\to\infty} \frac{f(x)}{\varphi(x)} = k_f \text{ exists} \Big\},
\]
with norm
\[
\|f\|_\varphi = \sup_{x\ge0} \frac{|f(x)|}{\varphi(x)}.
\]
\begin{theorem}[Korovkin-type convergence in $C_\varphi([0,\infty))$]
	\label{thm:korovkin-polynomial}
	Let \(n>\beta\) and \(\alpha>-1\). Then for every
	\(f\in C_\varphi([0,\infty))\),
	\[
	\lim_{n\to\infty}\|M_n^{(\alpha,\beta)}f-f\|_\varphi=0.
	\]
\end{theorem}

\begin{proof}
	We verify the Korovkin test functions \(1,t,t^2\). We note that
	\[
	M_n^{(\alpha,\beta)}1(x)=\sum_{k=0}^\infty\psi_{n,k}(x)=1\qquad(\forall x\ge0).
	\]
	Thus \(\|M_n^{(\alpha,\beta)}1-1\|_\varphi=0\).
	Also, 
	\[
	M_n^{(\alpha,\beta)} t(x)=\frac{n x + (\alpha+1)}{\,n-\beta\,}\,.
	\]
	Therefore
	\[
	M_n^{(\alpha,\beta)} t(x) - x = \frac{\beta x + \alpha+1}{n-\beta},
	\]
	and dividing by \(1+x^2\) gives
	\[
	\frac{|M_n^{(\alpha,\beta)} t(x)-x|}{1+x^2} \le \frac{1}{n-\beta}\sup_{y\ge0}\frac{|\beta y+\alpha+1|}{1+y^2}
	= \frac{C_1(\alpha,\beta)}{n-\beta},
	\]
	where \(C_1(\alpha,\beta)<\infty\). Thus \(\|M_n^{(\alpha,\beta)} t-t\|_\varphi\to0\).\\
	Recall 
	\[
	M_n^{(\alpha,\beta)} t^2(x)=\frac{n^2x^2 + n x(2\alpha+4) + (\alpha+1)(\alpha+2)}{(n-\beta)^2}.
	\]
	Compute the difference
	\[
	\begin{aligned}
		M_n^{(\alpha,\beta)} t^2(x)-x^2
		&= x^2\Big(\frac{n^2}{(n-\beta)^2}-1\Big) + \frac{n x(2\alpha+4)}{(n-\beta)^2} + \frac{(\alpha+1)(\alpha+2)}{(n-\beta)^2}.
	\end{aligned}
	\]
	Thus there exists a finite constant \(C_2(\alpha,\beta)\) with
	\[
	\frac{|M_n^{(\alpha,\beta)} t^2(x)-x^2|}{1+x^2}\le \frac{C_2(\alpha,\beta)}{n-\beta},
	\]
	for all \(x\ge0\). Hence \(\|M_n^{(\alpha,\beta)} t^2 - t^2\|_\varphi\to0\).\\
	By the weighted Korovkin theorem
	(positivity, linearity, and approximation of the three test functions),
	\(\|M_n^{(\alpha,\beta)} f - f\|_\varphi\to0\) for every \(f\in C_\varphi([0,\infty))\).
\end{proof}

\section{Local and Global $L_p$-Convergence Analysis }
In this section, we examine the convergence properties of the Szász–Mirakyan–Durrmeyer operators in the $L_p$-metric, both locally and globally on the half-line $[0,\infty)$. The analysis focuses on establishing sufficient conditions under which the operators provide approximation in $L_p$-spaces for $1\leq p < \infty$. Local convergence is studied on compact intervals, while global convergence is obtained through appropriate weight functions ensuring integrability over the unbounded domain. These results extend the classical $L_p$-approximation theory to the present class of operators and demonstrate their robustness in handling functions of varying growth behavior.
	\begin{theorem}[Korovkin-type theorem for $M_n^{(\alpha,\beta)}$]\label{thm:korovkin}
	For all $f \in C(K)$,
		\[
		\lim_{n \to \infty} \| M_n^{(\alpha,\beta)}[f] - f \|_{C(K)} = 0,
		\]
		i.e., $M_n^{(\alpha,\beta)}[f] \to f$ uniformly on every compact subset $K$ of $[0,\infty)$.
	\end{theorem}
	
	\begin{proof}
		We follow the classical Korovkin approach by verifying the convergence on the test functions $\{1, t, t^2\}$.
		From the moment calculations, we have
		\[
			M_n^{(\alpha,\beta)}[1](x) = 1.
			\]
		As $n \to \infty$, we have uniformly on compact sets
			\[
			M_n^{(\alpha,\beta)}[t](x)=\frac{nx + \alpha + 1}{n - \beta}  \to x.
			\]
			and 
			\[
			M_n^{(\alpha,\beta)}[t^2](x) =\frac{n^2x^2 + 2n(\alpha+1)x + (\alpha+1)(\alpha+2)}{(n-\beta)^2} \to x^2.
			\]
					Since \(M_n^{(\alpha,\beta)}\) is a sequence of positive linear operators and the convergence holds for \(1, t, t^2\), the classical Korovkin theorem yields
				\[
				\lim_{n \to \infty} \| M_n^{(\alpha,\beta_n)}[f] - f \|_{C(K)} = 0
				\]
		
		This completes the proof.
	\end{proof}
	
	
	\begin{corollary}
		Under the same conditions as Theorem \ref{thm:korovkin}, the operators $M_n^{(\alpha,\beta)}$ approximate continuous functions with exponential weights. That is, for any $f \in C([0,\infty))$ satisfying a growth condition $|f(t)| \leq Me^{At}$ for some $M,A > 0$, we have uniform convergence on compact subsets of $[0,\infty)$.
	\end{corollary}

In the study of approximation by positive linear operators on unbounded intervals, convergence may be investigated in different functional settings. In particular, for the operators \(M_n^{(\alpha,\beta)}\), it is natural to distinguish between two complementary notions of convergence in the \(L_p\)-sense. The first concerns the \emph{local} (or compact) \(L_p\)-convergence, which ensures that the approximation holds uniformly over every finite subinterval of \([0,\infty)\). The second involves the global \(L_p\)-convergence in weighted spaces equipped with an exponentially decaying weight, which controls the behavior of functions and operators at infinity. The results presented below address both these aspects, providing a comprehensive analysis of the approximation properties of the Sz\'{a}sz--Mirakyan--Laguerre--Durrmeyer operators.
\begin{lemma}\label{lem:E_n_bound}
	Let $\beta\ge0$ be fixed and let $\alpha\in[-\tfrac12,0]$. For $n>\beta$ and $t\ge0$ define
	\[
	E_n(t):=\frac{1}{n}(n-\beta)^{\alpha+1}t^{\alpha}\frac{\gamma(\alpha+1,(n-\beta)t)}{\Gamma(\alpha+1)},
	\]
	where $\gamma(s,z)=\int_0^z u^{s-1}e^{-u}\,du$ is the lower incomplete Gamma function.  
	Then for every fixed $R>0$ there exists a finite constant $C=C(R,\alpha,\beta)$ such that
	\[
	\sup_{n>\beta}\sup_{t\in[0,R]} E_n(t)\le C.
	\]
\end{lemma}

\begin{proof}
	Set $z=(n-\beta)t$ and define
	\[
	g(z):=\frac{z^{\alpha}\,\gamma(\alpha+1,z)}{\Gamma(\alpha+1)},\qquad z\ge0.
	\]
	A simple rearrangement gives
	\[
	E_n(t)=\Big(1-\frac{\beta}{n}\Big)g(z).
	\]
	Since $0<1-\beta/n\le1$ for all $n>\beta$, it suffices to show $\sup_{z\ge0} g(z)<\infty$.
	
	First consider the behavior as $z\to0^+$. Put $s=\alpha+1\in(\tfrac12,1]$. The small-argument expansion
	\[
	\gamma(s,z)=\frac{z^{s}}{s}+o(z^{s})\qquad(z\to0^+)
	\]
	yields
	\[
	g(z)=\frac{z^{\alpha}\gamma(\alpha+1,z)}{\Gamma(\alpha+1)}
	= \frac{z^{2\alpha+1}}{(\alpha+1)\Gamma(\alpha+1)}+o(z^{2\alpha+1}).
	\]
	For $\alpha\ge-1/2$ we have $2\alpha+1\ge0$, so $g(z)$ is finite at $z=0$ (indeed $g(0)=0$ when $\alpha>-1/2$, and $g(0)$ is finite when $\alpha=-1/2$).
	
	Next consider the behavior as $z\to\infty$. Since $\gamma(\alpha+1,z)\to\Gamma(\alpha+1)$, we have
	\[
	g(z)\sim z^{\alpha}\qquad(z\to\infty).
	\]
	Because $\alpha\le0$ the right-hand side is bounded as $z\to\infty$ (it tends to $0$ if $\alpha<0$ and to $1$ if $\alpha=0$). Therefore $g$ is bounded at infinity.
	
	On $(0,\infty)$ the function $g$ is continuous, hence the finiteness at both endpoints implies $g$ is bounded on $[0,\infty)$. Set
	\[
	C_0:=\sup_{z\ge0} g(z) <\infty.
	\]
	Then for every $n>\beta$ and every $t\in[0,R]$,
	\[
	E_n(t)=\Big(1-\frac{\beta}{n}\Big)g(z)\le g(z)\le C_0,
	\]
	so the claimed uniform bound holds with $C(R,\alpha,\beta):=C_0$.
\end{proof}
\begin{theorem}[Local $L_p$--convergence]
	\label{thm:local-Lp-fixed-beta}
	Let $1 \le p < \infty$, $-\frac{1}{2}\leq \alpha\leq 0$, and let $\beta \ge 0$ be a fixed real number such that $n > \beta$. Then for every fixed $R > 0$ and every $f \in L_p([0,R])$,
	\[
	\lim_{n \to \infty} \big\| M_n^{(\alpha,\beta)}[f] - f \big\|_{L_p([0,R])} = 0.
	\]
\end{theorem}

\begin{proof}
	Each $M_n^{(\alpha,\beta)}$ is a positive linear operator. By interchanging sum and integral (justified by Fubini's theorem and non-negativity), we obtain the kernel representation
	\[
	M_n^{(\alpha,\beta)}[f](x) = \int_0^\infty f(t) K_n(x,t)  dt,
	\]
	where the kernel is given by
	\[
	K_n(x,t) = e^{-nx} \sum_{k=0}^{\infty} \frac{(n - \beta)^{k + \alpha + 1}}{\Gamma(k + \alpha + 1)} t^{k+\alpha} e^{-(n - \beta)t} \frac{(nx)^k}{k!}.
	\]
	The kernel is non-negative for $x, t \ge 0$. Moreover, for $f \equiv 1$, we have
	\[
	M_n^{(\alpha,\beta)}[1](x) = 1 \quad \text{for all } x \ge 0,
	\]
	which implies that for each fixed $x \ge 0$,
	\[
	\int_0^\infty K_n(x,t)  dt = 1.
	\]
	Thus, $K_n(x, \cdot)$ is a probability density on $[0, \infty)$.\\
	Using theorem \eqref{thm:korovkin}, for every $g \in C([0,R])$,
	\[
	\lim_{n \to \infty} \| M_n^{(\alpha,\beta)}[g] - g \|_{\infty,[0,R]} = 0.
	\]
	Let $f \in L_p([0,R])$ and $\varepsilon > 0$. Since $C([0,R])$ is dense in $L_p([0,R])$, there exists $g \in C([0,R])$ such that
	\[
	\| f - g \|_{L_p([0,R])} < \varepsilon.
	\]
	
	We decompose the error as:
	\[
	\| M_n^{(\alpha,\beta)}[f] - f \|_{L_p([0,R])} \le \| M_n^{(\alpha,\beta)}[f - g] \|_{L_p([0,R])} + \| M_n^{(\alpha,\beta)}[g] - g \|_{L_p([0,R])} + \| g - f \|_{L_p([0,R])}.
	\]
	
	The third term is less than $\varepsilon$ by construction. The second term satisfies:
	\[
	\| M_n^{(\alpha,\beta)}[g] - g \|_{L_p([0,R])} \le R^{1/p} \| M_n^{(\alpha,\beta)}[g] - g \|_{\infty,[0,R]} \to 0 \quad \text{as } n \to \infty,
	\]
	since uniform convergence implies convergence in $L_p$ on bounded intervals.
	
	It remains to bound the first term. We claim there exists a constant $C = C(p,R,\alpha,\beta) > 0$, independent of $n$ for large $n$, such that for all $h \in L_p([0,R])$,
	\[
	\| M_n^{(\alpha,\beta)}[h] \|_{L_p([0,R])} \le C \| h \|_{L_p([0,R])}.
	\]
	
	To prove this, we use the Schur test \cite[Theorem 3.6]{zhu2007operator} for integral operators . Consider the kernel $K_n(x,t)$ restricted to $x, t \in [0,R]$. We verify:
	\begin{enumerate}
		\item $\sup_{x \in [0,R]} \int_0^R K_n(x,t)  dt \le 1$ (since $\int_0^\infty K_n(x,t)  dt = 1$ and $K_n \ge 0$).
		\item $\sup_{t \in [0,R]} \int_0^R K_n(x,t)  dx \le C$ for some constant $C$ independent of $n$.
	\end{enumerate}
	
	The first condition is immediate. For the second condition, we analyze
	\[
	I_n(t) = \int_0^R K_n(x,t)  dx = \int_0^R e^{-nx} \sum_{k=0}^{\infty} \frac{(n - \beta)^{k + \alpha + 1}}{\Gamma(k + \alpha + 1)} t^{k+\alpha} e^{-(n - \beta)t} \frac{(nx)^k}{k!}  dx.
	\]
	Interchanging sum and integral (justified by Tonelli's theorem):
	\[
	I_n(t) = \sum_{k=0}^{\infty} \frac{(n - \beta)^{k + \alpha + 1}}{\Gamma(k + \alpha + 1)} t^{k+\alpha} e^{-(n - \beta)t} \frac{n^k}{k!} \int_0^R x^k e^{-nx}  dx.
	\]
	Using the bound $\int_0^R x^k e^{-nx}  dx \le \int_0^\infty x^k e^{-nx}  dx = \frac{k!}{n^{k+1}}$, we obtain:
	\[
	I_n(t) \le \sum_{k=0}^{\infty} \frac{(n - \beta)^{k + \alpha + 1}}{\Gamma(k + \alpha + 1)} t^{k+\alpha} e^{-(n - \beta)t} \frac{1}{n^{k+1}} n^k = \frac{1}{n} \sum_{k=0}^{\infty} \frac{(n - \beta)^{k + \alpha + 1}}{\Gamma(k + \alpha + 1)} t^{k+\alpha} e^{-(n - \beta)t}.
	\]
	
	Now, observe that the series can be expressed using the confluent hypergeometric function:
	\[
	\sum_{k=0}^{\infty} \frac{(n - \beta)^{k + \alpha + 1}}{\Gamma(k + \alpha + 1)} t^{k+\alpha} e^{-(n - \beta)t} = (n - \beta)^{\alpha+1} t^\alpha e^{-(n - \beta)t} \sum_{k=0}^{\infty} \frac{[(n - \beta)t]^k}{\Gamma(k + \alpha + 1)}.
	\]
	Using the identity $\sum_{k=0}^{\infty} \frac{z^k}{\Gamma(k + \alpha + 1)} = e^z \frac{\gamma(\alpha+1, z)}{\Gamma(\alpha+1)}$, where $\gamma(\alpha+1, z)$ is the lower incomplete Gamma function, we get:
	\[
	I_n(t) \le \frac{1}{n} (n - \beta)^{\alpha+1} t^\alpha e^{-(n - \beta)t} \cdot e^{(n - \beta)t} \frac{\gamma(\alpha+1, (n - \beta)t)}{\Gamma(\alpha+1)} = \frac{1}{n} (n - \beta)^{\alpha+1} t^\alpha \frac{\gamma(\alpha+1, (n - \beta)t)}{\Gamma(\alpha+1)}.
	\]
	Using Lemma \ref{lem:E_n_bound},  there exists a constant $C = C(R,\alpha,\beta)$ such that
	\[
	\sup_{t \in [0,R]} I_n(t) \le C \quad \text{for all } n > \beta.
	\]
	Thus, by the Schur test \cite[Theorem 3.6]{zhu2007operator}, the operators $M_n^{(\alpha,\beta)}$ are uniformly bounded on $L_p([0,R])$:
	\[
	\| M_n^{(\alpha,\beta)}[h] \|_{L_p([0,R])} \le C \| h \|_{L_p([0,R])}.
	\]
	Returning to the error decomposition, we have for large $n$:
	\[
	\| M_n^{(\alpha,\beta)}[f] - f \|_{L_p([0,R])} \le C \varepsilon + o(1) + \varepsilon,
	\]
	where $o(1) \to 0$ as $n \to \infty$. Since $\varepsilon > 0$ is arbitrary, we conclude that
	\[
	\lim_{n \to \infty} \| M_n^{(\alpha,\beta)}[f] - f \|_{L_p([0,R])} = 0.
	\]
\end{proof}
We now introduce a weighted $L_p$--space suitable for analyzing the global behavior of the operators $M_n^{(\alpha,\beta)}$.

\begin{definition}[Weighted $L_p$--space]\label{def:weighted-Lp}
	Let $1 \le p < \infty$ and $\gamma \ge 0$.  
	The \emph{weighted space} $L_p^\gamma([0,\infty))$ consists of all measurable functions $f : [0,\infty) \to \mathbb{R}$ such that
	\[
	\|f\|_{L_p^\gamma}^p := \int_0^\infty |f(x)|^p e^{\gamma x}\,dx < \infty.
	\]
	The factor $e^{\gamma x}$ serves as an exponential weight, emphasizing the behavior of $f$ on the unbounded interval $[0,\infty)$.  
	For $\gamma = 0$, the space $L_p^\gamma([0,\infty))$ coincides with the usual $L_p([0,\infty))$ space.
\end{definition}
\begin{theorem}[Global $L_p$-convergence in weighted spaces]\label{thm:global-weighted-Lp-corrected}
	Let $1\le p<\infty$, let $\beta\ge0$ be fixed, and assume $\alpha\in [-\tfrac12,0].$	Fix $\gamma\ge0$ satisfying the condition 
	\[\gamma \le p\beta. 	\tag{H$_\gamma$}\]
	(If $\beta=0$ this forces $\gamma=0$, i.e.\ no exponential weight.) 
	Then for all sufficiently large $n>\beta$ the operators $M_n^{(\alpha,\beta)}$ map $L_p^\gamma([0,\infty))$ into itself with a uniform bound
	\[
	\|M_n^{(\alpha,\beta)}\|_{L_p^\gamma\to L_p^\gamma}\le C,
	\]
	where $C=C(\alpha,\beta,p,\gamma)$ is independent of $n$ (for all large $n$). Moreover for every $f\in L_p^\gamma([0,\infty))$,
	\[
	\lim_{n\to\infty}\|M_n^{(\alpha,\beta)}[f]-f\|_{L_p^\gamma}=0.
	\]
\end{theorem}

\begin{proof}
	Write the kernel representation
	\[
	M_n^{(\alpha,\beta)}[f](x)=\int_0^\infty K_n(x,t)\,f(t)\,dt,
	\]
	with
	\[
	K_n(x,t)=e^{-nx}\sum_{k=0}^\infty\frac{(n-\beta)^{k+\alpha+1}}{\Gamma(k+\alpha+1)}
	t^{k+\alpha}e^{-(n-\beta)t}\frac{(nx)^k}{k!},
	\]
	a nonnegative kernel and \(\int_0^\infty K_n(x,t)\,dt=1\) for each fixed \(x\).
	
	We conjugate by the weight as usual. For \(h(x)=e^{\gamma x/p}f(x)\) define
	\[
	\widetilde M_n[h](x):=e^{\gamma x/p}M_n\big[e^{-\gamma\cdot/p}h(\cdot)\big](x)
	=\int_0^\infty \widetilde K_n(x,t)\,h(t)\,dt,
	\]
	with
	\[
	\widetilde K_n(x,t):=e^{\gamma x/p}K_n(x,t)e^{-\gamma t/p}.
	\]
	Then \(\|M_n[f]\|_{L_p^\gamma}=\|\widetilde M_n[h]\|_{L_p}\). Thus it suffices to show \(\widetilde M_n\) is uniformly bounded on ordinary \(L_p([0,\infty))\) and that \(\widetilde M_n[h]\to h\) in \(L_p\) for \(h\) in the dense subspace corresponding to continuous compactly supported \(f\).
	
	We apply the Schur test. It is enough to find constants \(A,B\) independent of \(n\) (for all large \(n\)) such that
	\[
	\sup_{x\ge0}\int_0^\infty \widetilde K_n(x,t)\,dt \le A,\qquad
	\sup_{t\ge0}\int_0^\infty \widetilde K_n(x,t)\,dx \le B.
	\]
	Then \(\|\widetilde M_n\|_{L_p\to L_p}\le A^{1/q}B^{1/p}\) with \(1/p+1/q=1\).
	
	\textbf{Estimate of the first Schur integral (choice of \(A\)).} \\
	For fixed \(x\),
	\[
	\begin{aligned}
		\int_0^\infty \widetilde K_n(x,t)\,dt
		&= e^{\gamma x/p}\int_0^\infty K_n(x,t)e^{-\gamma t/p}\,dt\\
		&= e^{\gamma x/p}e^{-nx}\sum_{k=0}^\infty\frac{(n-\beta)^{k+\alpha+1}}{\Gamma(k+\alpha+1)}\frac{(nx)^k}{k!}
		\int_0^\infty t^{k+\alpha}e^{-(n-\beta+\gamma/p)t}\,dt.
	\end{aligned}
	\]
	Evaluating the inner Gamma-integral,
	\[
	\int_0^\infty t^{k+\alpha}e^{-(n-\beta+\gamma/p)t}\,dt
	=\frac{\Gamma(k+\alpha+1)}{(n-\beta+\gamma/p)^{k+\alpha+1}},
	\]
	and thus (after cancellations)
	\[
	\int_0^\infty \widetilde K_n(x,t)\,dt
	= \Big(\frac{n-\beta}{\,n-\beta+\gamma/p\,}\Big)^{\!\alpha+1}
	e^{\gamma x/p}e^{-nx}\sum_{k=0}^\infty\frac{(nx)^k}{k!}r_n^k,
	\]
	where \(r_n:=\dfrac{n-\beta}{\,n-\beta+\gamma/p\,}\in(0,1]\). Summing the exponential series gives
	\[
	\int_0^\infty \widetilde K_n(x,t)\,dt
	= \Big(\frac{n-\beta}{\,n-\beta+\gamma/p\,}\Big)^{\!\alpha+1}
	\exp\!\Big(\frac{\gamma x}{p}+nx(r_n-1)\Big).
	\]
	Note that
	\[
	r_n-1=-\frac{\gamma/p}{\,n-\beta+\gamma/p\,},
	\]
	so the exponent simplifies to
	\[
	\frac{\gamma x}{p}+nx(r_n-1)=x\cdot\frac{\gamma}{p}\cdot\frac{-\beta+\gamma/p}{\,n-\beta+\gamma/p\,}.
	\]
	Under hypothesis (H$_\gamma$), \(\gamma\le p\beta\), the factor \(-\beta+\gamma/p\le0\), hence the exponential factor is bounded by \(1\) for all \(x\ge0\). Also \(\big(\tfrac{n-\beta}{n-\beta+\gamma/p}\big)^{\alpha+1}\le1\) because \(\alpha+1>0\) and the denominator is at least the numerator. Consequently
	\[
	\int_0^\infty \widetilde K_n(x,t)\,dt \le 1\qquad\text{for all }x\ge0,\; n>\beta,
	\]
	and we may take \(A=1\).
	
	
	\textbf{Estimate of the second Schur integral (choice of \(B\)).} \\
	We must bound
	\[
	\int_0^\infty \widetilde K_n(x,t)\,dx
	= e^{-\gamma t/p}\int_0^\infty e^{\gamma x/p}K_n(x,t)\,dx.
	\]
	Interchange sum and integral (Tonelli) and compute the \(x\)-integral termwise. Using the identity (put \(u=nx\))
	\[
	\int_0^\infty e^{\gamma x/p}\frac{(nx)^k}{k!}e^{-nx}\,dx
	=\frac{1}{n}\int_0^\infty e^{-(1-\gamma/(np))u}\frac{u^k}{k!}\,du
	=\frac{1}{n}\big(1-\tfrac{\gamma}{np}\big)^{-(k+1)},
	\]
	which is finite for all sufficiently large \(n\) provided \(\gamma<np\). (For fixed \(\gamma\) this holds for all large \(n\).) Hence
	\[
	\begin{aligned}
		\int_0^\infty \widetilde K_n(x,t)\,dx
		&= e^{-\gamma t/p}\sum_{k=0}^\infty\frac{(n-\beta)^{k+\alpha+1}}{\Gamma(k+\alpha+1)}t^{k+\alpha}e^{-(n-\beta)t}\cdot\frac{1}{n}\big(1-\tfrac{\gamma}{np}\big)^{-(k+1)}\\[4pt]
		&= \frac{1}{n}(n-\beta)^{\alpha+1}t^\alpha e^{-(n-\beta)t}\big(1-\tfrac{\gamma}{np}\big)^{-1}
		\sum_{k=0}^\infty\frac{\big((n-\beta)t(1-\tfrac{\gamma}{np})^{-1}\big)^k}{\Gamma(k+\alpha+1)}.
	\end{aligned}
	\]
	As in the local analysis the series is identified with the incomplete Gamma factor and we obtain
	\[
	\int_0^\infty \widetilde K_n(x,t)\,dx
	\le \big(1-\tfrac{\gamma}{np}\big)^{-1}\;E_n\!\big(t(1-\tfrac{\gamma}{np})^{-1}\big),
	\]
	where
	\[
	E_n(u)=\frac{1}{n}(n-\beta)^{\alpha+1}u^{\alpha}\frac{\gamma(\alpha+1,(n-\beta)u)}{\Gamma(\alpha+1)}.
	\]
	By Lemma \ref{lem:E_n_bound} (which requires \(\alpha\in[-\tfrac12,0]\)) the function \(u\mapsto E_n(u)\) is uniformly bounded in \(n\) and \(u\ge0\). Therefore there exists \(C_1=C_1(\alpha,\beta)>0\) such that for all sufficiently large \(n\)
	\[
	\sup_{t\ge0}\int_0^\infty \widetilde K_n(x,t)\,dx \le \big(1-\tfrac{\gamma}{np}\big)^{-1}C_1 \le C,
	\]
	with \(C\) independent of \(n\) (for large \(n\)). Thus we may take \(B=C\).
	
	\textbf{Conclusion of boundedness and convergence.} \\
	With \(A=1\) and the above \(B\) we obtain a uniform Schur bound \(\|\widetilde M_n\|_{L_p\to L_p}\le B^{1/p}\), hence
	\[
	\|M_n\|_{L_p^\gamma\to L_p^\gamma}\le B^{1/p},
	\]
	for all sufficiently large \(n\). This proves the uniform operator bound.
	
	To prove convergence, fix \(f\in L_p^\gamma\). For \(R>0\) write \(f=f_1+f_2\) with \(f_1=f\chi_{[0,R]}\) and \(f_2=f\chi_{(R,\infty)}\). Given \(\varepsilon>0\) choose \(R\) so large that \(\|f_2\|_{L_p^\gamma}<\varepsilon\). Then by the uniform operator bound \(\|M_n f_2\|_{L_p^\gamma}\le B^{1/p}\varepsilon\). On the compact piece \(f_1\in L_p([0,R])\) and Theorem \ref{thm:local-Lp-fixed-beta} gives \(\|M_n f_1-f_1\|_{L_p^\gamma}\to0\) as \(n\to\infty\) (weights are bounded on \([0,R]\)). Combining these yields
	\[
	\limsup_{n\to\infty}\|M_n f-f\|_{L_p^\gamma}\le (1+B^{1/p})\varepsilon,
	\]
	and since \(\varepsilon>0\) is arbitrary the theorem follows.
\end{proof}
\section{Eigenfunctions} 
In this section, we study the spectral properties of the operators \(M_n^{(\alpha,\beta)}\), focusing on the characterization of their eigenvalues and associated eigenfunctions..
\begin{theorem}[Eigenpairs of \(M_n^{(\alpha,\beta)}\) and the coefficient matrix \(P\)]
	\label{thm:eigenpairs-Mn}
	Let \(n>\beta\) and \(\alpha>-1\). Define the operator $M_n^{(\alpha,\beta)}$ as defined in \ref{SMoperators}. Define the infinite matrix \(P=(P_{k,j})_{k,j\ge0}\) by
	\[
	P_{k,j}
	=\frac{(n-\beta)^{k+\alpha+1}n^j}{\Gamma(k+\alpha+1)j!}\,
	\frac{\Gamma(j+k+\alpha+1)}{(2n-\beta)^{\,j+k+\alpha+1}}
	=\Big(\frac{n-\beta}{2n-\beta}\Big)^{k+\alpha+1}\Big(\frac{n}{2n-\beta}\Big)^j
	\frac{(k+\alpha+1)_j}{j!}.
	\]
	Then \(P\) is nonnegative and row-stochastic, and the operator \(M_n^{(\alpha,\beta)}\) and \(P\) are related as follows: if \(v=(v_j)_{j\ge0}\) is any coefficient vector for which the Poisson series
	\[
	\Phi_v(x):=\sum_{j\ge0} v_j\,\psi_{n,j}(x)
	\]
	converges, then
	\[
	M_n^{(\alpha,\beta)}[\Phi_v]=\Phi_{Pv}.
	\]
	Consequently the following eigenpairs hold for \(M_n^{(\alpha,\beta)}\):
	\begin{enumerate}
		\item \( \lambda_1=1\) is an eigenvalue with eigenfunction \(\phi_1(x)\equiv1\).
		\item \( \lambda_2=\big(1-\tfrac{\beta}{n}\big)^{\alpha+1}\) is an eigenvalue with eigenfunction \(\phi_2(x)=e^{-\beta x}\).
	\end{enumerate}
	More precisely, the vectors \(v^{(1)}=(1,1,1,\dots)\) and \(v^{(2)}=\big(1-\tfrac{\beta}{n}\big)^j_{j\ge0}\) satisfy
	\(Pv^{(1)}=v^{(1)}\) and \(Pv^{(2)}=\lambda_2 v^{(2)}\), and lifting these vectors via \(\Phi\) yields the stated eigenfunctions of \(M_n^{(\alpha,\beta)}\).
\end{theorem}

\begin{proof}  First, we prove that, one can interchange the sum and the integral in
	\[
	\int_0^\infty\Big(\sum_{j\ge0} v_j\psi_{n,j}(t)\Big)\,t^{k+\alpha}e^{-(n-\beta)t}\,dt
	=\sum_{j\ge0} v_j\int_0^\infty \psi_{n,j}(t)\,t^{k+\alpha}e^{-(n-\beta)t}\,dt,
	\]
	using the Fubini--Tonelli theorem (Folland \cite[Theorem 2.18]{folland1999real} or Royden--Fitzpatrick \cite[Theorem 17.7]{royden1988real}). Two convenient sufficient hypotheses are:
	
	\begin{enumerate}
		\item (Nonnegative coefficients) If $v_j\ge0$ for all $j$, then the integrand $\sum_{j\ge0} v_j\psi_{n,j}(t)\,t^{k+\alpha}e^{-(n-\beta)t}$ is nonnegative and Tonelli's theorem applies, permitting the interchange.
		\item (Bounded coefficients) If $\sup_{j\ge0}|v_j|=:C<\infty$, then for every $t\ge0$
		\[
		\Big|\sum_{j\ge0} v_j\psi_{n,j}(t)\Big|\le C\sum_{j\ge0}\psi_{n,j}(t)=C.
		\]
		Hence
		\[
		\Big|\sum_{j\ge0} v_j\psi_{n,j}(t)\,t^{k+\alpha}e^{-(n-\beta)t}\Big|
		\le C\,t^{k+\alpha}e^{-(n-\beta)t},
		\]
		and the right-hand side is integrable on $[0,\infty)$ because $\int_0^\infty t^{k+\alpha}e^{-(n-\beta)t}\,dt<\infty$ (recall $n>\beta$, $\alpha>-1$). Thus by the dominated convergence theorem (or Fubini's theorem for absolutely integrable integrands) we may swap the sum and integral.
	\end{enumerate}
	
	
	 Expand an arbitrary function \(\Phi_v(x)=\sum_{j\ge0}v_j\psi_{n,j}(x)\) (with coefficients such that the series converges). Using above argument and the definition of \(M_n^{(\alpha,\beta)}\) we compute
	\[
	\begin{aligned}
		M_n^{(\alpha,\beta)}[\Phi_v](x)
		&=\sum_{k\ge0}\left(\frac{(n-\beta)^{k+\alpha+1}}{\Gamma(k+\alpha+1)}
		\int_0^\infty\Big(\sum_{j\ge0} v_j\psi_{n,j}(t)\Big)\,t^{k+\alpha}e^{-(n-\beta)t}\,dt\right)\psi_{n,k}(x)\\
		&=\sum_{k\ge0}\Big(\sum_{j\ge0}P_{k,j}v_j\Big)\psi_{n,k}(x)
		=\Phi_{Pv}(x),
	\end{aligned}
	\]
	where
	\[
	P_{k,j}=\frac{(n-\beta)^{k+\alpha+1}}{\Gamma(k+\alpha+1)}\cdot\frac{n^j}{j!}
	\int_0^\infty t^{j+k+\alpha} e^{-(2n-\beta)t}\,dt,
	\]
	and evaluation of the Gamma integral yields the above formula for \(P_{k,j}\). Summing the series for \(\sum_j P_{k,j}\) (or invoking \(\sum_j\psi_{n,j}(t)=1\)) shows each row sums to \(1\), hence \(P\) is row-stochastic and nonnegative.\\
	For the constant eigenfunction, take \(v^{(1)}=(1,1,\dots)\). Since each row of \(P\) sums to \(1\) we have \(Pv^{(1)}=v^{(1)}\). Lifting gives
	\(\Phi_{v^{(1)}}(x)=\sum_{j\ge0}\psi_{n,j}(x)=1\) and therefore \(M_n^{(\alpha,\beta)}[1]=1\).
	
	For the second eigenpair consider the geometric vector \(v^{(2)}\) with \(v^{(2)}_j=z^j\) and seek \(z\) such that \(Pv^{(2)}=\lambda v^{(2)}\). A direct summation (using the Pochhammer/binomial series)
	\[
	\sum_{j\ge0}\frac{(k+\alpha+1)_j}{j!}a^j=(1-a)^{-(k+\alpha+1)},\qquad |a|<1,
	\]
	with \(a=\dfrac{n z}{2n-\beta}\), yields
	\[
	(Pv^{(2)})_k=\Big(\frac{n-\beta}{2n-\beta-n z}\Big)^{k+\alpha+1}.
	\]
	Equating this to \(\lambda z^k\) for all \(k\) forces \(z\) to satisfy
	\[
	z=\frac{n-\beta}{2n-\beta-n z},
	\]
	whose solutions are \(z=1\) and \(z=1-\dfrac{\beta}{n}\). The \(z=1\) case reproduces the constant eigenvector. For \(z=1-\beta/n\) the resulting eigenvalue is
	\[
	\lambda_2=\Big(\frac{n-\beta}{2n-\beta-n z}\Big)^{\alpha+1}
	=\Big(\frac{n-\beta}{n}\Big)^{\alpha+1}=\Big(1-\frac{\beta}{n}\Big)^{\alpha+1}.
	\]
	Lifting \(v^{(2)}\) with \(z=1-\beta/n\) gives
	\[
	\Phi_{v^{(2)}}(x)=\sum_{j\ge0}\Big(1-\frac{\beta}{n}\Big)^j\psi_{n,j}(x)
	=e^{-n x}\sum_{j\ge0}\frac{\big(n x (1-\beta/n)\big)^j}{j!}=e^{-\beta x},
	\]
	so \(e^{-\beta x}\) is an eigenfunction of \(M_n^{(\alpha,\beta)}\) with eigenvalue \(\lambda_2\). This completes the proof.
\end{proof}
\begin{remark}
	The operator \( M_n^{(\alpha,\beta)} \) is positive and linear. The constant eigenfunction corresponds to the invariant property of \( M_n^{(\alpha,\beta)} \). The second eigenvalue
	\(\lambda_2 = \left(1 - \frac{\beta}{n}\right)^{\alpha + 1}\)
	satisfies \( 0 < \lambda_2 < 1 \) for \( 0 < \beta < n \),
	and therefore the component along \( e^{-\beta x} \) decays geometrically under iteration:
	\[
	(M_n^{(\alpha,\beta)})^{r}[e^{-\beta x}]
	= \lambda_2^{\,r} e^{-\beta x}, \qquad r \ge 1.
	\]
\end{remark}
\end{document}